\documentclass[11pt]{article}
\usepackage{authblk} 

\usepackage[latin1]{inputenc}
\usepackage[T1]{fontenc}
\usepackage{amsmath}
\usepackage{amsfonts}
\usepackage{amssymb}
\usepackage{enumerate}

\usepackage[normalem]{ulem} 

\usepackage{amsmath,xcolor,ulem}
\usepackage{amsfonts}
\usepackage{amssymb,todonotes}
\usepackage{amsthm}
\usepackage{graphicx}
\usepackage{epsfig}
\usepackage{marginnote}

\usepackage[colorlinks=true, allcolors=blue]{hyperref}

\usepackage[top=2.8cm,bottom=2.8cm,left=2.5cm,right=2.5cm]{geometry}
\usepackage{float}
\usepackage[algoruled]{algorithm2e}
\usepackage[font=footnotesize,width=\linewidth]{caption}

\newtheorem{thm}{Theorem}

\newcommand{\C} {\mathbb{C}}
\newcommand{\R} {\mathbb{R}}

\newcommand{\hu} {\hat{u}}

\newcommand{\iu}{\mathrm{i}} 
\renewcommand{\Re}{\operatorname{Re}\,} 
\renewcommand{\Im}{\operatorname{Im}\,} 

\DeclareMathOperator{\g}{g}

\title{Absorbing Boundary Conditions for Variable Potential Schr\"odinger Equations via Titchmarsh-Weyl Theory}

\author{Matthias Ehrhardt\footnote{Corresponding author, \href{mailto:ehrhardt@uni-wuppertal.de}{ehrhardt@uni-wuppertal.de}} ,
Chunxiong Zheng\footnote{Supported by the Alexander-von-Humboldt Foundation,  \href{mailto:czheng@tsinghua.edu.cn}{czheng@tsinghua.edu.cn}} 
}

\affil{IMACM, School of Mathematics and Natural Sciences, \\ University of Wuppertal, Germany}
\affil{Department of Mathematical Sciences, \\ Tsinghua University, Beijing 100084, P.R. China}

\begin{document}
\maketitle

\begin{tikzpicture}[remember picture,overlay]
	\node[anchor=north east,inner sep=20pt] at (current page.north east)
	{\includegraphics[scale=0.2]{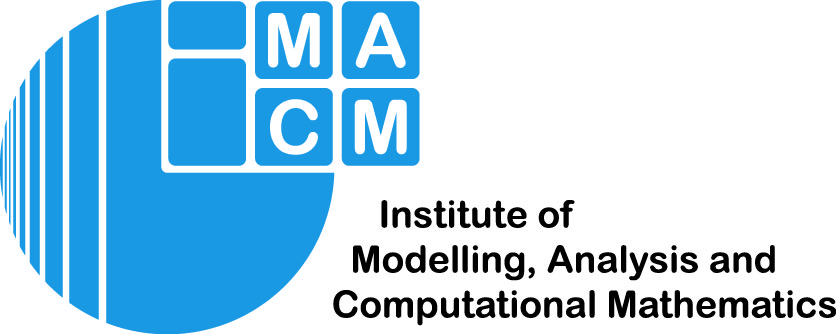}};
\end{tikzpicture}

\begin{abstract}
We propose a novel approach to simulate the solution of the time-dependent Schr\"odinger equation with a general variable potential. 
The key idea is to approximate the Titchmarsh-Weyl m-function (exact Dirichlet-to-Neumann operator) by a rational function with respect to an appropriate spectral parameter. 
By using this method, we overcome the usual high-frequency restriction associated with absorbing boundary conditions in general variable potential problems. 
The resulting fast computational algorithm for absorbing boundary conditions ensures accuracy over the entire frequency range.
\end{abstract}

\begin{minipage}{0.9\linewidth}

\medskip

\end{minipage}

\section{Introduction}
In this paper we consider the linear Schr\"odinger problem of the form
\begin{equation}\label{sch}
\begin{split}
   \iu u_t+\partial_x^2u&=V(x)u,\quad (x,t)\in\R\times(0,T],\\
   u(x,0)&=u_0(x),\quad x\in\R,
\end{split}
\end{equation}
where $T$ denotes the finite evolution time, and $u_0$ is an initial wave packet supported in a finite interval $\Omega_{\rm int}=[x_-,x_+]$ with $x_-<x_+$. 
It is well known that under mild conditions the Cauchy problem \eqref{sch} has a unique solution
$u\in C(\R^{+},L^2(\R))$, cf.\ \cite{Pazy}, e.g.:
\begin{thm}
Let $u_0\in L^2(\R)$ and real-valued potential $V\in L^\infty(\R)$.
Then the problem \eqref{sch} has a unique solution $u\in C(\R^{+},L^2(\R))$.
Moreover, the ``energy'' is preserved, i.e.\
\begin{equation}
   \|u(.,t)\|_{L^2(\R)}=\|u_0\|_{L^2(\R)},\quad\forall\,t \ge 0.
\end{equation}
\end{thm}

The Schr\"odinger problem \eqref{sch} is defined on an unbounded domain $x\in\R$.
To numerically simulate its solution, it is common practice to truncate the domain to a bounded one,
for example, $\Omega_{{\rm int},T}=\Omega_{\rm int}\times(0,T]$.
\emph{Absorbing boundary conditions} (ABCs) are thus necessary for well-posedness at the two artificially introduced boundaries, $\Sigma_{\pm,T}=\{x_\pm\}\times(0,T]$.

Numerical simulation of the linear Schr\"odinger equation on
unbounded domains with an external potential has been a hot
research area for nearly thirty years, cf.\ the concise review article \cite{TBC}.
An ABC is called \emph{exact} if the solution of the truncated domain problem remains the same
as that of the original unbounded domain problem.
The exact ABC is guaranteed to exist due to the well-posedness
of the linear Schr\"odinger problem \eqref{sch},
but it can only be formulated analytically for some special potentials,
such as constant potential \cite{EhAr01}, linear potential \cite{EhMi04}, symmetric periodic potential \cite{EZ},
isotropic free particle potential,
Morse potential, harmonic potential, and Bargeman potential,
cf.\ e.g.\ \cite{Pol02}.
In the more general case, i.e., for general variable potential problems
one is led to design approximate analytical ABCs for a given frequency regime with respect to some a priori criterion.
Methods in this category include the pseudo-differential calculus method \cite{AnBeMo04,AnBeDe06, AnBeSz09}, the perfectly matched layer (PML) method \cite{ZhengPML07}, 
and the operator splitting method \cite{ZXW08}.
To the authors' knowledge, all of them are essentially based on the \emph{high frequency approximations}.
For low-frequency problems, the ABCs would be less accurate by these methods.

This paper proposes a new approach to the design of ABCs for the Schr\"odinger problem.
Inspired by the work of Alpert, Greengard, and Hag\-strom \cite{AlpertGreengardHagstrom} on the fast evaluation of nonreflecting boundary kernels for time-domain wave propagation,
we approximate the \emph{Titchmarsh-Weyl m-function} (equivalently, the exact DtN operator) in the frequency domain by a rational function with respect to an appropriate spectral parameter.
In the time domain, the nonreflecting boundary kernels are thus approximated by a sum of exponentials, which makes the approximate ABCs easy to implement.

The rationality of the above treatment is due to the analyticity property and the asymptotic behavior of the m-function.
Since our approximation is performed in the whole frequency regime, the proposed ABCs are expected to be more versatile and accurate, 
especially in the low-frequency regime, thus overcoming the typical high-frequency restriction.
Note that the Titchmarsh-Weyl m-function is nothing else but the so-called \emph{total symbol} in microdifferential calculus,
which is treated by an asymptotic expansion to obtain a hierarchy of ABCs, cf.\ \cite{AnBeMo04,AnBeDe06, AnBeSz09}.
Note also that Titchmarsh-Weyl theory is used in the analysis of initial value problems for Schr\"odinger equations operator-valued potentials \cite{GZ13} or strongly singular potentials \cite{AK12}. It is also used in practical applications in quantum mechanics \cite{Chiba90, GDB98} and in option pricing in mathematical finance \cite{LiZhang}.

This paper is organized as follows.
 In Section~\ref{s:2} we review the basic facts of the Titchmarsh-Weyl theory for Schr\"odinger operators in one dimension for ease of later reference.
Then, in Section~\ref{s:3} we discuss the Titchmarsh-Weyl m-function (i.e.\ the exact Dirichlet-to-Neumann operator) and explain the algorithm used to compute the m-function numerically.
Thus, at least from a numerical point of view, the exact ABC is explicitly known, see Section~\ref{s:3}. 
However, when simulating the Schr\"odinger equation \eqref{sch}, the difficulty does not lie in the computation of the m-function in a frequency domain method, 
which is presented in Section~\ref{s:4}, but in its inverse Laplace transformation, which is too expensive.
For this reason, in Section~\ref{s:5} we introduce a rational approximation of the m-function in the frequency domain
to obtain an approximate ABC that can be computed efficiently using a fast evaluation technique \cite{ZhengFast} in the time domain.
We discuss practical implementation issues and finally in Section~\ref{s:6} we conclude with numerical results
illustrating that our new approach leads to an efficient and reliable algorithm for the time-dependent Schr\"odinger equation with a general variable potential.

\section{The Titchmarsh-Weyl theory}\label{s:2}
We will review here, for convenience of later work, the essentials of
Titchmarsh-Weyl (TW) theory for Schr\"odinger operators in one dimension.
The interested reader may consult \cite[Section~2]{GZ} or \cite[Section~4.3]{Benn20} for a more detailed presentation.

First, consider the Schr\"odinger operator $L$ on the real line given by
\begin{equation}\label{L}
    L=-\partial_x^2+V(x),\qquad x\in\R,
\end{equation}
with a real-valued, locally integrable potential $V\in L^1_{\rm loc}(\R)$.
Let $x_0\in\R$ be an arbitrarily chosen point, called \emph{reference point}.
In the following we will study how the solutions depend on this parameter $x_0$.

To do this, we consider $\theta(x;x_0,\lambda)$ and $\phi(x;x_0,\lambda)$ to be the \emph{fundamental solutions} of the Schr\"odinger eigenvalue problem.
\begin{equation}\label{schrod}
     -u_{xx}+V(x)u=\lambda u,\qquad x\in\R,\quad\lambda\in\C,
\end{equation}
with the following initial conditions at the reference point $x_0$:
 \begin{subequations}
\begin{align}
\theta(x_0;x_0,\lambda)&=1,\qquad \theta_x(x_0;x_0,\lambda)=0,\\
  \phi(x_0;x_0,\lambda)&=0,\qquad \phi_x(x_0;x_0,\lambda)=1.
\end{align}
\end{subequations}
It can be shown that under these assumptions
$\theta(x;x_0,\lambda)$ and $\phi(x;x_0,\lambda)$ exist on the whole real axis,
and they are entire functions of $\lambda$ and real for $\lambda\in\R$.
Now, as a basic fact of TW theory, the equation \eqref{schrod} has at least one solution $\psi_{\pm}$, called \emph{Weyl's solution} with
\begin{subequations}\label{psicond}
\begin{equation}\label{psicond1}
    \psi_{\pm}(x_0;x_0,\lambda)=1,
\end{equation}
and
\begin{equation}\label{psicond2}
   \psi_{\pm}(x;x_0,\lambda)\in L^2(\R^{x_0}_{\pm})
\end{equation}
for any $\lambda\in\C_+$.
\end{subequations}
Here $\R^{x_0}_{\pm}$ stands for the interval $[x_0,\pm\infty)$ and $\C_+$ denotes the upper half of the complex plane, i.e.\ $\C_+=\{z\in\C\, |\Im z>0\}$.
A potential $V(x)$ is said to be in the {\em limit-point case} at $\pm\infty$ if and only if there exists only one Weyl's solution in the corresponding $L^2$ space.
The reader will immediately realize that the assumption of $V(x)$ in the limit-point case is necessary for the well-posedness of the Schr\"odinger problem \eqref{sch} in a more general setting.
At positive infinity, a standard sufficient condition for the limit-point case is given by \cite{RS2}:
\begin{thm}[{\cite[Theorem X.8]{RS2}}]
Let $V(x)$ be a continuous real-valued function on $(x_0,\infty)$ and suppose that there
exists a positive differentiable function $M(x)$ such that
\begin{enumerate}
  \item[(i)] $V(x)\ge-M(x)$ if $x>x_0$;
  \item[(ii)] $\int\limits_{x_1}^\infty\bigl(M(x)\bigr)^{-1/2}\,dx=\infty$ for any $x_1>x_0$;
  \item[(iii)] $M'(x)/\bigl(M(x)\bigr)^{3/2}$ is bounded near $\infty$.
\end{enumerate}
Then $V(x)$ is in the limit-point case at  $\infty$.
\end{thm}
An analogous result can be given at negative infinity point.

According to this theorem, a potential $V(x)$ is in the limit-point case
provided that $V(x)\ge-k x^2$ for some constant $k$ and for all sufficiently large $x$.
This implies that the restriction on the potential for the limit-point case is very weak:
 It only excludes some especially strange potentials, which may not be physically relevant at all.
Roughly speaking, the limit case does not admit potentials that tend too fast
(faster than quadratically) to $-\infty$ for $x\to\pm\infty$.

Due to the boundary conditions \eqref{psicond1} we can write
\begin{equation}\label{psim}
   \psi_{\pm}(x;x_0,\lambda)= \theta(x;x_0,\lambda)+m_{\pm}(x_0,\lambda)\,\phi(x;x_0,\lambda),
\end{equation}
with some uniquely determined coefficient,
the \emph{Titchmarsh-Weyl m-function} $m_{\pm}(x_0,\lambda)$.
This function plays a fundamental role in the spectral theory of
the Schr\"odinger operator \eqref{L} on the half-line $\R^{x_0}_{\pm}$.

We will now summarize some of the most important properties of the
Titchmarsh-Weyl m-function. First,
\begin{equation}\label{m1}
    m_{\pm}(x_0,\lambda) \text{ is analytic with respect to } \lambda \text{ to } \C\backslash\R\
\text{ and } m_{\pm}\colon\C_+\to\C_+
\end{equation}
and is therefore called a \emph{Herglotz function} (or Nevanlinna or Pick function),
cf.\ \cite[Lemma~2.3]{GZ}.
It is easy to show that this Herglotz-property is directly related to the positive type of the DtN-map in the sense of memory equations,
cf.\ \cite{EhAr01} for the corresponding case of constant external potential.
Thus it is an essential ingredient of the stability w.r.t.\ the $L^2$ norm.

We also have the \emph{symmetry property}
\begin{equation}\label{m2}
   \overline{m_{\pm}(x_0,\lambda)}=m_{\pm}(x_0,\bar{\lambda})
\end{equation}
and the local singularities of $m$
are real and and most of them are first order, i.e.\
\begin{align}
\lim_{\epsilon\to0+}(-\iu\epsilon)\,m_{\pm}(x_0,\lambda+\iu\epsilon)&\ge0,\quad\lambda\in\R,
\end{align}
cf.\ \cite[Theorem~A.2]{GZ}.

Another important property is given by the Borg-Marchenko theorem \cite{Borg, March}, 
which states that the Titchmarsh-Weyl m-function $m_{\pm}(x_0,\lambda)$ 
uniquely determines the potential $V(x)$ at $x>x_0$ (or $x<x_0$).
Moreover, since $\psi_\pm(x;x_0,\lambda)$ changes with a simple multiplication 
when the reference point $x_0$ changes, one has
\begin{equation}\label{psim3}
   m_{\pm}(x,\lambda)= \frac{\partial_x\psi_{\pm}(x;x_0,\lambda)}{\psi_{\pm}(x;x_0,\lambda)}.
\end{equation}
Thus, it is easy to verify that the $m$ function satisfies the following \emph{Riccati} equation:
\begin{equation}\label{Ricc}
    \partial_xm_{\pm}(x,\lambda)=-m_{\pm}^2(x,\lambda)+V(x)-\lambda.
\end{equation}

\section{The exact ABC by Titchmarsh-Weyl theory}\label{s:3}
We apply the Laplace transform
\begin{equation}\label{Laplacedef}
    \hu(x,s) = \mathcal{L}\bigl(u(x,t)\bigr)(s) = \int_0^{+\infty } u(x,t) \,e^{-st}\,dt, \quad \Re s>0,
\end{equation}
to the Schr\"odinger equation \eqref{sch} on the right exterior domain $\Omega_+=\{x\in\R|\,x>x_+\}$ and on the left exterior domain $\Omega_-=\{x\in\R|\,x<x_-\}$.
In the frequency domain, the Schr\"odinger equation is a second order homogeneous ODE
\begin{equation}\label{schLaplace}
    -\hu_{xx}+V(x)\hu=\lambda\hu,\quad x\in\Omega_\pm,
\end{equation}
with $\lambda=\iu s\in\C_+$. The exact absorbing boundary condition of the DtN form in the frequency domain is thus
\begin{equation*}
    \hat{u}_x(x_\pm,\lambda) = m_\pm(x_\pm,\lambda)\hat{u}(x_\pm,\lambda).
\end{equation*}

Only in some special cases, the m-function has a closed analytic form \cite{Brown,Simon}.
For example, in the case of a constant potential $V\equiv V_0$ one gets
\begin{equation}\label{mSE}
    m_+(x_+,\lambda)=-\sqrt[+]{-\lambda+V_0}.
\end{equation}
If the potential represents a \emph{harmonic oscillator}, i.e.\ $V(x)=x^2$ on the interval $[0,\infty)$,
one obtains a meromorphic m-function given by the ratio of two gamma functions:
\begin{equation}\label{mGamma}
m_+(0,\lambda)=-\frac{2\Gamma(\frac{3}{4}-\frac{1}{4}\lambda)}
{\Gamma(\frac{1}{4}-\frac{1}{4}\lambda)}.
\end{equation}
Finally, for the \emph{Bargmann potential}
\begin{equation}\label{VBargeman}
    V(x)=-8\beta^2\,\frac{\beta-\gamma}{\beta+\gamma}\frac{e^{-2\beta x}}{(1+\frac{\beta-\gamma}{\beta+\gamma}e^{-2\beta x})^2},\quad\beta>0,\quad\gamma\ge 0,
\end{equation}
one obtains the m-function
\begin{equation}\label{mBargeman}
    m_+(0,\lambda)=-\sqrt[+]{-\lambda}-\frac{\gamma^2-\beta^2}{\sqrt[+]{-\lambda}+\gamma}.
\end{equation}

In the general case, however, numerical methods must be considered.
This problem has been studied in many papers, e.g.\ \cite{Brown, WGK99, GKW00, Kirby}.
In this paper we simply compute the $m$-function by evolving the Riccati equation \eqref{Ricc} 
with the classical fourth-order Runge-Kutta scheme and setting an initial data $m_\pm(x_{\pm,\lambda},\lambda)=\mp\sqrt[+]{-\lambda}$ 
at a sufficiently distant point $x_{\pm,\lambda}=\pm 200$. 
This treatment is reasonable since the potentials in our numerical tests actually decay to zero for $x\to\infty$.

\section{The frequency-domain method}\label{s:4}
The solution of the time-dependent Schr\"odinger equation could then be computed 
using the following frequency-domain method:
\begin{enumerate}[Step 1.]
\item Fix $\sigma>0$.
For each $s=\sigma+\iu\mu$ with $\mu\in\R$, solve the
Laplace-transformed Schr\"odinger equation in the bounded interval 
$[x_-,x_+]$:
\begin{align*}
  -\hu_{xx}+V(x)\hu &= \iu s\,\hu -\iu u_0(x),\quad x\in[x_-,x_+],\\
  \hu_x(x_-)&=m_-(x_-,\iu s)\,\hu(x_-),\\
  \hu_x(x_+)&=m_-(x_+,\iu s)\,\hu(x_+).
\end{align*}
\item Perform the inverse Laplace transformation
\begin{equation}\label{step2}
\begin{split}
u(x,t) &= \mathcal{L}^{-1}\bigl(\hat{u}(s,t)\bigr)(x)\\
&=\frac{1}{2\pi
\iu}\int_{\sigma-\iu\infty}^{\sigma+\iu\infty}e^{st}\hu(s,t)\,ds
=
\frac{e^{\sigma t}}{2\pi}\int_{-\infty}^{\infty}e^{\iu ft}u(x,\sigma+\iu f)\,df,
\end{split}
\end{equation}
to derive the wave function $u(x,t)$ for any $t\in(0,T]$.
\end{enumerate}
In the numerical implementation, some parameters have to be tuned.
The function $\hat{u}$ is smoother for larger \emph{damping factor} $\sigma$, but the evolution time span is then limited
because an exponential term is involved in the inverse Laplace transformation. 
As a common practice, we set $\sigma=1/T$, where $T$ is a prescribed evolution time. 
The integration domain is unbounded in \eqref{step2} and must be truncated. 
We introduce a cut-off frequency $f_c$ and confine the integration to the interval $[-f_c,f_c]$. 
In addition, to get rid of high frequency oscillations in the inverse-transformed function, we should introduce another filtering function $\chi$,
which remains 1 over a sufficiently large frequency band with zero frequency as its center, and vanishes smoothly near the endpoints of $[-f_c,f_c]$. 
A good candidate (empirically) is
\begin{equation*}
    \chi=\exp\bigl(-(1.2f/f_c)^{20}\bigr).
\end{equation*}
After these treatments, we then derive an approximate inverse transformation as
\begin{equation}\label{cutoff}
  \frac{e^{\sigma t}}{2\pi}\int_{-\infty}^{\infty}e^{\iu ft}u(x,\sigma+\iu f)\,df\approx
  \frac{e^{\sigma t}}{2\pi}\int_{-f_{c}}^{f_{c}}\chi(f)\,e^{\iu ft}u(x,\sigma+\iu f)\,df.
\end{equation}
The right side is computed with an appropriate quadrature scheme.

\section{The time-domain method}\label{s:5}
It follows from Section~\ref{s:2} that the exact ABC we are looking for is now explicitly known, at least from a numerical point of view.
But this is not the whole story for simulating the solution of the time-dependent Schr\"odinger equation.
The difficulty lies not in computing the m-function itself, but in computing its inverse Laplace transformation.
Of course, a numerical inverse transformation is possible, but it would be too expensive.

Therefore, in this section we design an approximate ABC based on the rational approximation of the m-function.
The kernel functions are of exponential type with respect to the half-order time derivative operator, 
so the fast evaluation technique proposed in \cite{ZhengFast} (cf.\ Appendix) is applicable.
For some alternative fast evaluation methods, we refer the reader to \cite{TBC} and the references therein.
The rational approximation is realized by solving a least squares problem, 
an analogous technique to that used in \cite{AlpertGreengardHagstrom} for fast evaluation of the boundary kernel functions of the hyperbolic wave equation.

In the time domain, the truncated Schr\"odinger problem reads
\begin{equation}
\begin{split}
  \iu u_t+\partial_x^2u &= V(x)u,\quad (x,t)\in[x_-,x_+]\times(0,T],\\
  u(x,0)&=u_0(x),\quad x\in[x_-,x_+],\\
  u_x(x_\pm,t) &= \mathcal{L}^{-1}\bigl(m_\pm(x_\pm,\iu s)\,\hat{u}(x_\pm,s)\bigr)(t),\quad t\in (0,T].
\end{split}
\end{equation}
To simplify the notation, we will focus on the right boundary at $x=x_+$.
Let us recall that the DtN map in the frequency domain is
\begin{equation*}
   \hat{u}_x(x_+,s) = m_+(x_+,\iu s)\,\hat{u}(x_+,s).
\end{equation*}
Returning to the time domain we have to consider the convolution
\begin{equation*}
  u_x(x_+,t) = \mathrm{K}(t)*u(x_+,t),
  \quad\text{with}\quad \mathrm{K}(t) = \mathcal{L}^{-1}\bigl(m_+(x_+,\iu s)\bigr)(t).
\end{equation*}
There are two major difficulties here.
First, it is generally \emph{hard to compute} $\mathrm{K}(t)$
and second, the convolution involved naturally leads to a \emph{nonlocal-in-time DtN map}.

To get an idea, let us first consider two specific simple examples.
In the case of the free Schr\"odinger ($V\equiv 0$), cf.\
\eqref{mSE}, we have
\begin{equation}
  m_+(\iu s)=-\sqrt[+]{-\iu s}\quad\text{ and thus }\quad
\mathcal{L}^{-1}\bigl(m_+(\iu s)\bigr)=-e^{-\iu\pi/4}\,\partial_t^{\frac{1}{2}},
\end{equation}
with the half-order time derivative defined as
\begin{equation}\label{halftime}
    \partial_t^{\frac{1}{2}} v(t) = \frac{1}{\sqrt{\pi}}\, \frac{\rm d}{\rm d t}
    \int_0^t \frac{v(\tau)}{\sqrt{t-\tau}}\, d\tau,
\end{equation}
which can be efficiently evaluated with some existing methods,
e.g.\ \cite{ArEhSo03,JiGr04,ZhengFast}.
Second, we consider the Bargmann potential \eqref{VBargeman}, and in this case
the m-function is
\begin{equation*}
   m_+(0,is)=-\sqrt[+]{-\iu s}-\frac{\gamma^2-\beta^2}{\sqrt{-\iu s}+\gamma},    
\end{equation*}
thus we have
\begin{equation*}
    \mathcal{L}^{-1}\bigl(m_+(0,\iu s)\bigr) = -e^{-\iu\pi/4}\,\partial_t^{\frac{1}{2}}
    -(\gamma^2-\beta^2)(e^{-\iu\pi/4}\,\partial_t^{\frac{1}{2}}+\gamma)^{-1}.
\end{equation*}
This operator can then be evaluated efficiently by introducing an unknown function and using the fast methods for $\partial_t^{\frac{1}{2}}$.

Inspired by these two examples it naturally leads us to think about the possibility of
approximating the $m$-function with a rational function with respect to a new spectral parameter\ $k=\sqrt[+]{-\iu s}$ (NOT $s$), i.e.,
\begin{equation}\label{idea}
  m_+(x_+,\iu s)\approx \tilde{m}_+(x_+,\iu s)
  =-\sqrt[+]{-\iu s}+\sum_{n=1}^{d}\frac{\alpha_{n}}{\sqrt{-\iu s}+\beta_{n}}.
\end{equation}
Once this is done, we can then replace the exact $m$-function with the approximate alternative $\tilde{m}$,
which leads to the approximate kernel function
\begin{equation*}
   \mathcal{L}^{-1}\bigl(\tilde{m}_+(x_+,\iu s)\bigr) 
   = -e^{-\iu\pi/4}\,\partial_t^{\frac{1}{2}}+\sum_{n=1}^d\alpha_n(e^{-\iu\pi/4}\,\partial_t^{\frac{1}{2}}+\beta_n)^{-1}.
\end{equation*}
The analogous idea for the Bargmann potential \eqref{VBargeman} can then be used to handle this kernel function.

The answer for the possibility is affirmative considering that the \emph{asymptotic expansion} has been given in \cite{Danielyan} as
\begin{equation}\label{asym}
  m_+(x_+,\lambda)=-\sqrt[+]{-\lambda}+\mbox{o}(1/\sqrt{r}),\quad r\to\infty,
\end{equation}
where $\lambda=\mu r$, $r\in\R$ and the convergence is uniform for $\mu$ in any compact subset of $\C_+$, cf.\ \cite[Theorem~C.4]{Sims07}.

Now putting
\begin{equation}
  \g_+(\lambda)=m_+(x_+,\lambda)+\sqrt[+]{-\lambda},
\end{equation}
in view of \eqref{asym} we know that $\g_+(\lambda)$ is analytic in $\C_+$ with respect to $\lambda$ and
it tends to zero for $\lambda\to\infty$.
We then use the method of Alpert, Greengard and Hagstrom  \cite{AlpertGreengardHagstrom} to approximate $\g_+(\lambda)$ with a rational function with respect to $\sqrt[+]{-\lambda}$ (NOT $\lambda$).
In terms of \eqref{asym} we consider the following \emph{nonlinear least squares problem}
\begin{equation}\label{nonlin}
  \epsilon=\min_{P,Q}\int_{-\infty+\iu\sigma}^{\infty+\iu\sigma}
  \biggl|\frac{P(\sqrt[+]{-\lambda})}{Q(\sqrt[+]{-\lambda})}-\g_+(\lambda)\biggr|^2\,|d\sqrt[+]{-\lambda}|,
\end{equation}
where $P$, $Q$ are polynomials with $\deg(P)+1=\deg(Q)=d$, and $d$ is determined by making $\epsilon\le \epsilon_0$, where $\epsilon_0$ is a prescribed tolerance number.
This nonlinear problem \eqref{nonlin} is then solved using the technique of
 linearization and orthogonalization \cite{AlpertGreengardHagstrom}.
Finally, by expressing $P/Q$ as a sum of poles, we arrive at
\begin{equation}\label{approx}
m_+(x_+,\lambda)\approx \tilde{m}_+(x_+,\lambda)=-\sqrt[+]{-\lambda}+\sum_{n=1}^{d}\frac{\alpha_{n}}{\sqrt[+]{-\lambda}+\beta_{n}}.
\end{equation}
Note that the coefficients $\alpha_n$ and $\beta_n$ should appear as conjugate pairs due to the symmetry property \eqref{m2}.
Unfortunately, it is not clear whether the rational approximation $\tilde{m}_+(x_+,\lambda)$ in \eqref{approx} still has the important Herglotz-property of the m-function.
Moreover, the Herglotz-property cannot be checked by some conditions on the poles due to the leading square root in \eqref{approx}.

\begin{figure}[htbp]
\centering
\includegraphics[width=2.3in]{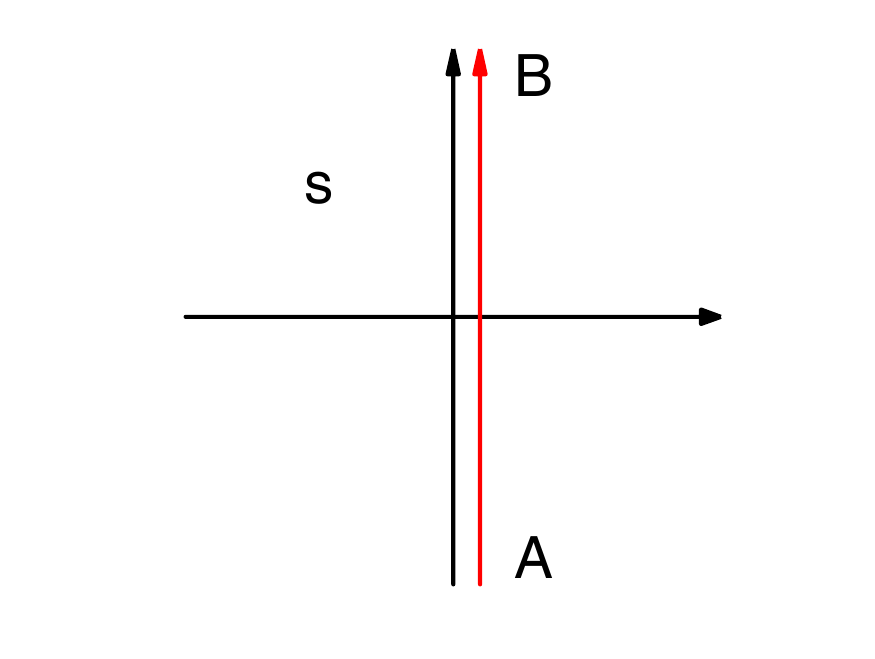}
\includegraphics[width=2.3in]{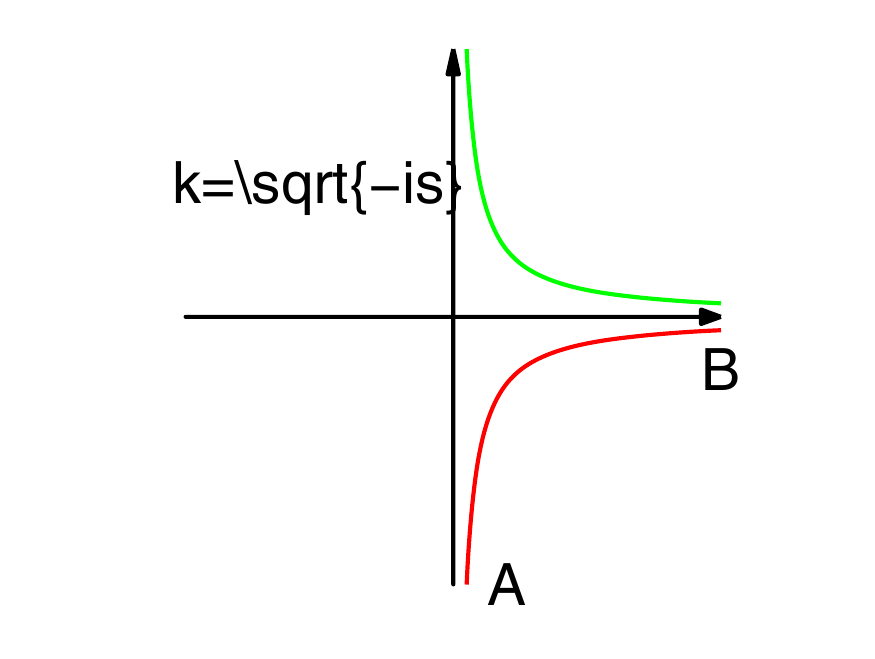}
\caption{$s$-plane and $k$-plane.}
\label{f1}
\end{figure}

Applying the same idea to the $m$-function $m_-(x_-,\lambda)$ and we get
the approximate boundary condition in the frequency domain
\begin{equation}\label{abcu}
    \hu_x(x_\pm,s) = \biggl(\mp \sqrt[+]{-\iu s}+\sum_{n=1}^{d_\pm}\frac{\alpha_{n,\pm}}{\sqrt[+]{-\iu s}+\beta_{n,\pm}}\biggr)\hu(x_\pm,s).
\end{equation}
If we introduce new unknowns $\hat{w}_{n,\pm}$ as
\begin{equation}
   \hat{w}_{n,\pm} = \frac{\hu(x_\pm,s)}{\sqrt[+]{-\iu s}+\beta_{n,\pm}},
\end{equation}
then we can rewrite \eqref{abcu} as
\begin{subequations}
\begin{align}
   \hu_x(x_\pm,s)\pm \sqrt[+]{-\iu s}\,\hu &= \sum_{n=1}^{d_\pm}\alpha_{n,\pm}\hat{w}_{n,\pm},\\
   \sqrt[+]{-\iu s}\,\hat{w}_{n,\pm}+\beta_{n,\pm}\hat{w}_{n,\pm} &= \hu(x_\pm,s),\qquad n=1,\dots,d_\pm.
\end{align}
\end{subequations}
In the time domain, the approximate boundary condition is
\begin{subequations}
\begin{align}
  u_x(x_\pm,t)\pm e^{-\iu\pi/4}\,\partial_t^{\frac{1}{2}}u(x_\pm,t)
  &= \sum_{n=1}^{d_\pm}\alpha_{n,\pm}w_{n,\pm}(t),\\
  e^{-\iu\pi/4}\,\partial_t^{\frac{1}{2}}w_{n,\pm}(t)+\beta_{n,\pm}w_{n,\pm}(t)&=u(x_\pm,t),\qquad n=1,\dots,d_\pm.
\end{align}
\end{subequations}
The final approximate truncated time-domain problem is formulated as
\begin{equation}
\begin{split}
  \iu u_t+\partial_x^2u&=V(x)u,\quad (x,t)\in\Omega_{\rm int}\times(0,T],\\
  u(x,0)&=u_0(x),\quad x\in\Omega_{\rm int},\\
  u_x(x_\pm,t)&\pm e^{-\iu\pi/4}\,\partial_t^{\frac{1}{2}}u(x_\pm,t)
  = \sum_{n=1}^{d_\pm}\alpha_{n,\pm}w_{n,\pm}(t),\\
  e^{-\iu\pi/4}\,\partial_t^{\frac{1}{2}}w_{n,\pm}(t)&+\beta_{n,\pm}w_{n,\pm}(t)
=u(x_\pm,t),\qquad n=1,\dots,d_\pm.
\end{split}
\end{equation}

\section{Numerical results}\label{s:6}
In this section, we present some numerical results to test the accuracy of the proposed methods.
In each example, the standard Crank-Nicolson scheme for time discretization is used.
The fast evaluation of the half-order time derivative operator \eqref{halftime} is performed with the method of Zheng \cite{ZhengFast}.
The computational domain is chosen to be $\Omega_{\rm int}=[x_-,x_+]=[-5.5]$ and the initial data is a Gaussian beam: $u_0(x)=e^{-x^2+4\iu x}$.
We use an $8^{\rm th}$-order FEM method with 1024 elements for the spatial discretization and a uniform time step of size $\Delta t=10^{-4}$.

\subsection{The Free Schr\"odinger Equation}
The exact solution of the free Schr\"odinger equation $(V(x)\equiv0)$ is
\begin{equation*}
   u_{\rm exa}(x,t) = \sqrt{\frac{\iu}{-4t+\iu}}\exp\biggl(\frac{-\iu x^2-4x+16t}{-4t+\iu}\biggr).
\end{equation*}
We set the cut-off frequency $f_c$ in \eqref{cutoff} to be $f_c=256$ and used a filtering function $\chi(f)=\exp\bigl(-(1.2f/f_{c})^{20}\bigr)$.
The following Table~\ref{t1} shows the relative $L^2$-errors to the exact solution at certain time points when using 8097 quadrature points with Simpson's rule.
\begin{table}[htbp]
\centering
\begin{tabular}{||c|c|c|c|c|c||}
\hline
Time points&0.5&0.6&0.7&0.8&0.9\\
\hline
Relative $L^2$-errors&2.26e-7&3.46e-8&7.60e-9&5.60e-9&6.25e-9\\
\hline
\end{tabular}
\caption{
Relative $L^2$-error at certain time points for the free Schr\"odinger equation.}
\label{t1}
\end{table}
Using this simple example, we can see that the frequency method with truncation and filtering works quite well: the magnitude of the relative $L^2$-errors is at most on the order of $10^{-7}$.

\subsection{The Coulomb-like Potential}
\begin{figure}[htbp]
\centering
\includegraphics[width=0.85\textwidth]{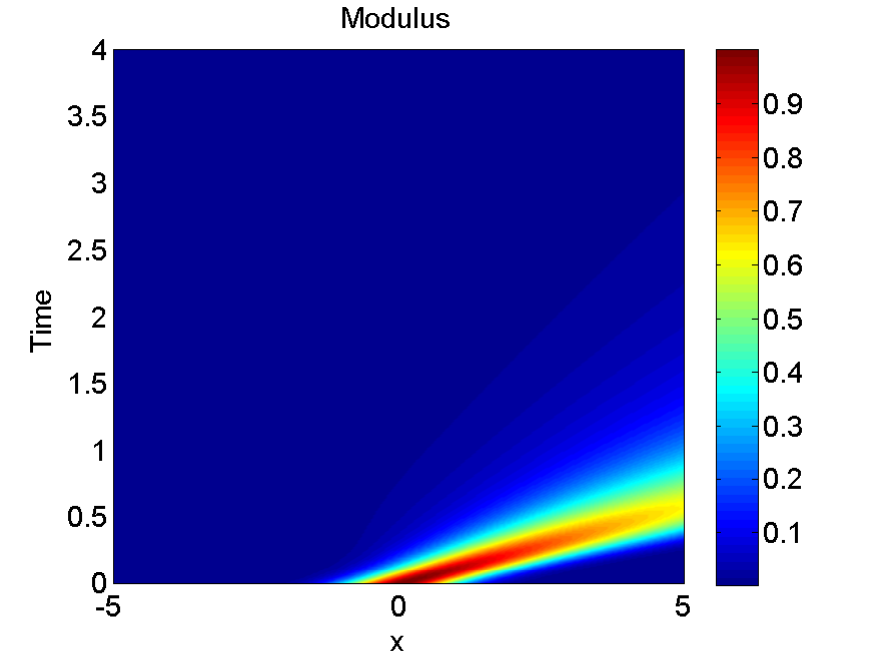}
\caption{Time evolution of solution with Coulomb potential \eqref{cp}.}
\label{c1}
\end{figure}

\begin{figure}[htbp]
\centering
\includegraphics[width=0.8\textwidth]{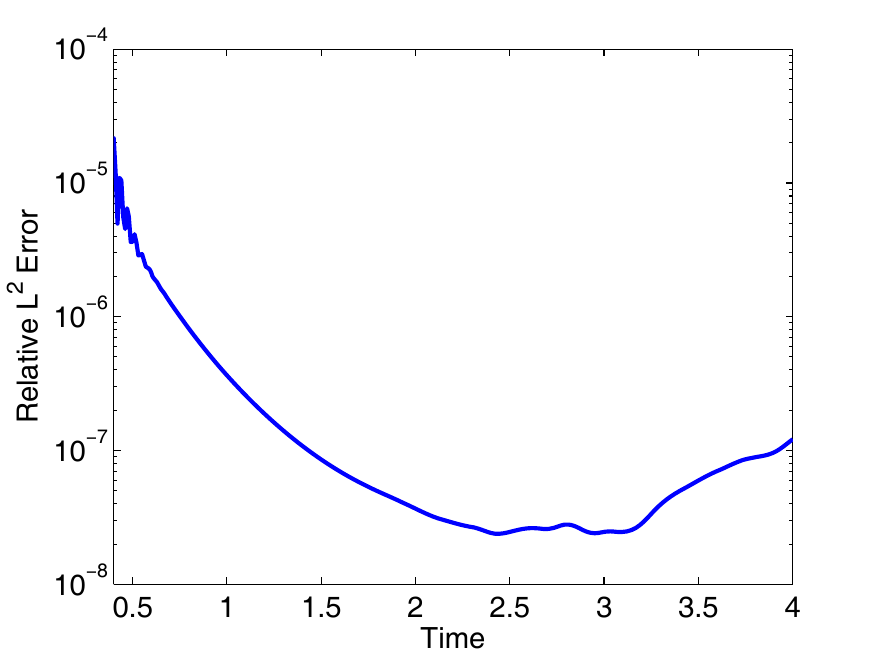}
\caption{Time evolution of the relative $L^2$-error.}
\label{cerror}
\end{figure}

In the second example we test the time-domain method with a Coulomb-like potential
\begin{equation}\label{cp}
   V(x)=\frac{1}{\sqrt{1+x^2}}.
\end{equation}
We set $\sigma$ for Step~1 to $\sigma=1$ and the tolerance number for the nonlinear least squares problem \eqref{nonlin} to $\epsilon_0=10^{-8}$. 
Here we get 4 poles.
Figure~\ref{c1} shows the time evolution in a colored contour plot. 
One can see how the initial beam spreads out as time increases.
Figure~\ref{cerror} shows how the relative $L^2$-error evolves in time.
In this example with a varying external potential the magnitude
of the relative $L^2$-errors remains below $10^{-5}$.

\subsection{The Gaussian Barrier}
\begin{figure}[htb]
\centering
\includegraphics[width=0.85\textwidth]{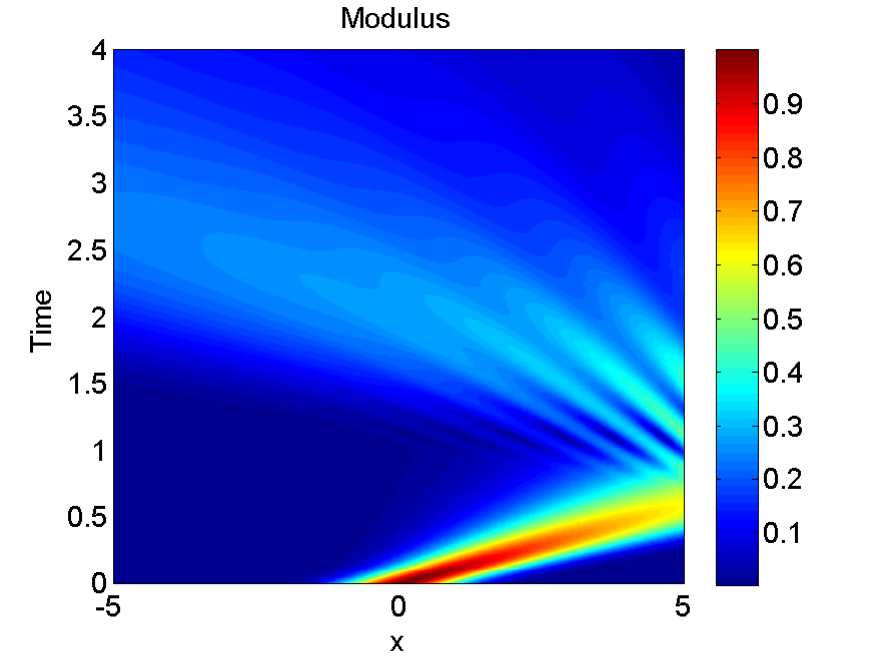}
\caption{Time evolution of solution with Gaussian barrier \eqref{gb}.}
\label{f:gb}
\end{figure}

\begin{figure}[htb]
\centering
\includegraphics[width=0.8\textwidth]{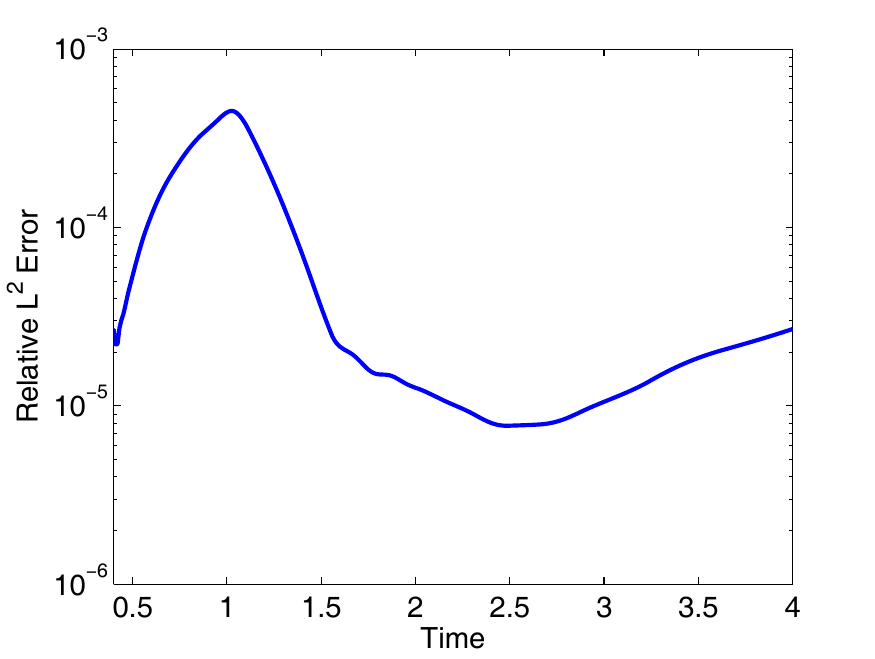}
\caption{Time evolution of the relative $L^2$-error.}
\label{gberror}
\end{figure}

Next, we change the potential to a Gaussian barrier
\begin{equation}\label{gb}
   V(x)=30e^{-36(x-8)^2}
\end{equation}
with a height of 30, located in the exterior domain $x>5$ and centered at $x=8$. 
We set $\sigma=1$ and $\epsilon_0=10^{-4}$ and get 21 poles with the nonlinear least squares algorithm. 
Figure~\ref{f:gb} shows the time evolution of the solution. 
One can clearly see how the initial beam propagates, spreads, and is (partially) reflected by the Gaussian barrier \eqref{gb}. 
The time evolution of the corresponding relative $L^2$-error is shown in Figure~\ref{gberror}.
The relative $L^2$ errors remain below $5\times 10^{-4}$.

Unfortunately, the nonlinear least squares algorithm used in this paper failed to produce a rational approximation within an error tolerance much smaller than $\epsilon_0$. 
A more efficient algorithm is still desired, and this problem is currently under investigation.

\section*{Conclusion and Outlook}
In this work we presented a new approach to simulate the solution to the Schr\"odinger equation with a general space-dependent potential in unbounded domains. 
Both frequency-domain and time-domain methods have been developed.

Future work will consist of implementing a more sophisticated algorithm for computing the $m$-function.
Instead of solving the Riccati equation \eqref{Ricc}, we will consider computing the Weyl circles \cite{Brown, Kirby}
or the boundary control approach of Avdonin, Mikhaylov and Rybkin \cite{AMR07}.
We will also look for a more stable algorithm for its rational approximation.
It will also be clarified how this rational approximation can be made to preserve the essential 
Herglotz-property of the analytic m-function.
This study will allow us to perform a rigorous stability analysis of this new approach.
Finally, as a future goal, we want to extend our approach to the multi-dimensional Schr\"odinger problem, following the idea of Amrein and Pearson \cite{AmPe04}.

\newpage
\section*{Appendix: Fast Evaluation Method}\label{A:ZhengFast}
Here we present a short description of the method in \cite{ZhengFast} for evaluating the
half-order time derivative $\partial_t^{\frac{1}{2}}$. 
For any smooth function $v=v(t)$ with $v(0)=v'(0)=0$, 
it is known that the semi-discrete half-order derivative
\begin{equation}\label{dishalf}
  D_t^{\frac{1}{2}}v(t_n):=\sqrt{\frac{2}{\Delta t}}\sum_{m=0}^n\alpha_m \,v(t_{n-m})
\end{equation}
with
\begin{equation}\label{ALPHA} 
\alpha_m=\begin{cases}
\beta_k=\frac{(2k)!}{2^{2k}(k!)^2}&,\ m=2k,\\
-\beta_k&,\ m=2k+1,
\end{cases}
\end{equation}
gives a second-order approximation of $\partial_t^{\frac{1}{2}}v(t_n)$ (see \cite{AntoineBesse2003,ZhengFast}).
Suppose there exists a sum of decaying exponentials satisfying
\begin{equation}\label{APP} 
\tilde{\beta}_k = \sum_{j=1}^Lw_j\,e^{-s_jk},\quad s_j>0,\
   |\beta_k-\tilde{\beta}_k|\le \epsilon,\ k=0,1,\dots,\bigg[\frac{N}{2}\bigg].
\end{equation}
Here $N$ is the total number of time steps. 
If $\epsilon$ is small enough, it is reasonable to approximate \eqref{dishalf} with
\begin{equation}\label{sum1}
   \tilde{D}_t^{\frac{1}{2}}v(t_n) := \sqrt{\frac{2}{\Delta t}}\bigl(v(t_n)-v(t_{n-1})\bigr) + \sqrt{\frac{2}{\Delta t}}\sum_{m=2}^n\tilde{\alpha}_mv(t_{n-m}),
\end{equation}
where
\begin{equation}\label{TALPHA} 
\tilde{\alpha}_m=
\begin{cases}
 \tilde{\beta}_k&,\ m=2k,\\
 -\tilde{\beta}_k&,\ m=2k+1.
\end{cases}
\end{equation}
Set $v_k=v(t_k)$, $\mathbf{v}=(v_0,v_1,\dots)$, and define
\begin{equation*}
\mathcal{F}_{\rm odd}(w,s;\mathbf{v},k)
:=\sum_{m=1}^kwe^{-sm}v_{2k+1-2m}
\end{equation*}
and
\begin{equation*}
\mathcal{F}_{\rm even}(w,s;\mathbf{v},k)
:=\sum_{m=1}^kwe^{-sm}v_{2k-2m}.
\end{equation*}
Thus,
$\mathcal{F}_{\rm odd}(w,s;\mathbf{v},0)=\mathcal{F}_{\rm even}(w,s;\mathbf{v},0)=0$.
In addition, we have the following recursions
\begin{eqnarray*}
\mathcal{F}_{\rm odd}(w,s;\mathbf{v},k)
&=&e^{-s} \bigl[wv_{2k-1}+\mathcal{F}_{\rm odd}(w,s;\mathbf{v},k-1)\bigr],\\
\mathcal{F}_{\rm even}(w,s;\mathbf{v},k)
&=&e^{-s}\bigl[wv_{2k-2}+\mathcal{F}_{\rm even}(w,s;\mathbf{v},k-1)\bigr].
\end{eqnarray*}
The summation \eqref{sum1} is then computed within $O(L)$ operations
as
\begin{equation*}
\sum_{m=2}^n\tilde{\alpha}_m v_{n-m}=
\begin{cases}
\displaystyle\sum_{j=1}^L\mathcal{F}_{\rm even}(w_j,s_j;\mathbf{v},k)-\sum_{j=1}^L\mathcal{F}_{\rm odd}(w_j,s_j;\mathbf{v},k-1)&, n=2k,\\
\displaystyle\sum_{j=1}^L\mathcal{F}_{\rm odd}(w_j,s_j;\mathbf{v},k)-\sum_{j=1}^L\mathcal{F}_{\rm even}(w_j,s_j;\mathbf{v},k)&,
n=2k+1.
\end{cases}
\end{equation*}
In \cite{ZhengFast} for $N=1,000,000$, the authors found a sum of 81
decaying exponentials that approximates $\beta_k$ with an error of less than $5.0\times 10^{-11}$.



\begin{thebibliography}{00}

\bibitem{akto}
T.~Aktosun,
\emph{Mathematical Studies in Nonlinear Wave Propagation},
Contemporary Mathematics \textbf{379},
 American Mathematical Society, Providence, RI, 2005.

\bibitem{AlpertGreengardHagstrom}
B.~Alpert, L.~Greengard and T.~Hagstrom,
\emph{Rapid evaluation of nonreflecting boundary kernels for time-domain wave propagation},
SIAM J. Numer. Anal. \textbf{37} (2000), 1138--1164.

\bibitem{AmPe04}
W.O.~Amrein and D.B.~Pearson,
\emph{M operators: a generalisation of Weyl-Titchmarsh theory},
J. Comput. Appl. Math. \textbf{171} (2004), 1--26.


\bibitem{AntoineBesse2003}
X. Antoine and C. Besse, Unconditionally stable discretization
schemes of non-reflecting boundary conditions for the
one-dimensional Schr\"odinger equation,
 J. Comput. Phys. \textbf{181} (2003), 157--175.

\bibitem{AnBeMo04}
X.~Antoine, C.~Besse and V.~Mouysset,
\emph{Numerical Schemes for the simulation of the two-dimensional
Schr\"odinger equation using non-reflecting boundary conditions},
Math.\ Comp.\ \textbf{73} (2004), 1779--1799.

\bibitem{AnBeDe06}
X.~Antoine, C.~Besse and S.~Descombes,
\emph{Artificial boundary conditions for one-dimensional cubic nonlinear Schr\"odinger equations},
SIAM J.\ Numer.\ Anal.\ \textbf{43} (2006), 2272--2293.

\bibitem{TBC}
X.~Antoine, A.~Arnold, C.~Besse, M.~Ehrhardt and A.~Sch\"adle,
\emph{A Review of Transparent and Artificial Boundary Conditions Techniques for Linear and Nonlinear Schr\"odinger Equations},
Commun. Comput. Phys. \textbf{4} (2008), 729--796. 

\bibitem{AnBeSz09}
X.~Antoine, C.~Besse and J.~Szeftel,
\emph{Towards accurate artificial boundary conditions for nonlinear PDEs through examples},
CUBO \textbf{11} (2009), 29--48.

\bibitem{ArEhSo03}
A.~Arnold, M.~Ehrhardt and I.~Sofronov,
\emph{Discrete transparent boundary conditions for the {S}chr{\"o}dinger
equation: fast calculation, approximation, and stability},
Comm.\ Math.\ Sci.\ \textbf{1} (2003), 501--556.




\bibitem{AMR07}
S.~Avdonin, V.~Mikhaylov and A.~Rybkin,
\emph{The boundary control approach to the Titchmarsh-Weyl m-function}, 
Comm. Math. Phys. \textbf{275} (2007),  791--803.

\bibitem{Benn80}
C.~Bennewitz and W.N. Everitt,
\emph{Some remarks on the Titchmarsh-Weyl m-coefficient.}
Proc. Pleijel Conference, Uppsala, Sweden, 1979, pp. 49--108.


\bibitem{Benn20}
C.~Bennewitz, R.~Weikard and M.~Brown,
\textit{Spectral and scattering theory for ordinary differential equations}, Springer, 2020.


\bibitem{Borg}
G.~Borg,
\emph{Uniqueness theorems in the spectral theory of $y''+(\lambda-q(x))y=0$},
Proc. 11th Scandinavian Congress of Mathematicians, Johan Grundt Tanums Forlag, Oslo,
1952, pp. 276--287.

\bibitem{Brown}
B.M.~Brown, V.G.~Kirby and J.D.~Pryce,
\emph{A numerical method for the determination of the Titchmarsh-Weyl m-coefficient},
Proc. R. Soc. Lond. A  \textbf{435} (1991),  535--549.





\bibitem{Chiba90}
Y.~Chiba and S.~Ohnishi,
\emph{Resonance-state calculation applying the Weyl-Titchmarsh theory: 
Application for the quantum-confined Stark effects on excitons in a GaAs--Al$_x$Ga$_{1-x}$As quantum well},
Phys. Rev. B \textbf{41} (1990), 6065--6068



\bibitem{Codd55}
E.A.~Coddington and N. Levinson,
\emph{Theory of Ordinary Differential Equations},
McGraw-Hill, New York, 1955.


\bibitem{Danielyan}
A.A.~Danielyan and B.M.~Levitan,
\emph{On the asymptotic behavior of the Weyl-Titchmarsh m-function},
Math. USSR Izvestiya  \textbf{36} (1991), 487--496.


\bibitem{EhAr01}
M.~Ehrhardt and A.~Arnold,
\emph{Discrete transparent boundary conditions for the Schr\"odinger equation},
Riv.\ Mat.\ Univ.\ Parma  \textbf{6} (2001), 57--108.

\bibitem{EhMi04}
M.~Ehrhardt and R.E.~Mickens,
\emph{Solutions to the discrete Airy equation: Application to parabolic equation calculations},
J.~Comput.\ Appl.\ Math.\ \textbf{172} (2004), 183--206.

\bibitem{Ehr07}
M.~Ehrhardt,
\emph{Discrete transparent boundary conditions for Schr\"odinger-type equations
for non-compactly supported initial data},
Appl.\ Numer.\ Math.\ \textbf{58} (2008), 660--673.

\bibitem{EZ}
M.~Ehrhardt and C.~Zheng,
\emph{Exact artificial boundary conditions for problems with periodic structures},
J.\ Comput.\ Phys.\ \textbf{227} (2008), 6877--6894.

\bibitem{Zheng}
M.~Ehrhardt, H.~Han and C.~Zheng,
\emph{Numerical simulation of waves in periodic structures},
 Commun.\ Comput.\ Phys.\ \textbf{5} (2009), 849--870.



\bibitem{Simon}
F.~Gesztesy and B.~Simon,
\emph{A new approach to inverse spectral theory, II. General real potentials
and the connection to the spectral measure},
Ann. Math.\  \textbf{152} (2000),  593--643.

\bibitem{GZ}
F.~Gesztesy and M.~Zinchenko,
\emph{On spectral theory for Schr\"odinger operators with strongly singular potentials},
Math. Nachr. \textbf{279} (2006), 1041--1082.

\bibitem{GZ13}
F.~Gesztesy, R.~Weikard and M.~Zinchenko,
\textit{Initial value problems and Weyl-Titchmarsh theory for Schr\"odinger operators with operator-valued potentials},
Oper. Matrices \textbf{7(2)} (2013), 241--283.



\bibitem{GDB98}
G.N.~Gibson, G.~Dunne and K.J.~Bergquist,
\emph{Tunneling Ionization Rates from Arbitrary Potential Wells},
Phys. Rev. Lett. \textbf{81} (1998), 2663--2666.

\bibitem{GKW00}
N.~Gordon, J.P.~Killingbeck and M.R.M.~Witwit,
\emph{Numerical Determination of the Titchmarsh-Weyl m-coefficient},
Comput. Math. \& Math. Phys. \textbf{40}, (2000), 108--117.



\bibitem{JiGr04}
S.~Jiang and L.~Greengard,
\emph{Fast evaluation of nonreflecting boundary conditions for the Schr\"odinger
equation in one dimension},
Comp.\ Math.\ Appl.\ \textbf{47} (2004), 955--966.


\bibitem{Kirby}
V.G.~Kirby,
\emph{A numerical method for determining the Titchmarsh-Weyl m-coefficient and its application
to certain integro-differential inequalities},
Ph.D. Thesis, University of Wales, College of Cardiff, UK, 1990.

\bibitem{AK12}
A.~Kostenko, A.~Sakhnovich and G.~Teschl, 
\textit{Weyl-Titchmarsh theory for Schr\"odinger operators with strongly singular potentials}, 
Int. Math. Res. Not. \textbf{8} (2012), 1699--1747.


\bibitem{LiZhang}
Y.~Li and J.E~Zhang,
\emph{Option pricing with Weyl-Titchmarsh theory},
Quant. Finance \textbf{4} (2004), 457--464.

\bibitem{March}
V.A.~Marchenko,
\emph{Certain problems in the theory of second-order differential operators},
Doklady Akad. Nauk SSSR \textbf{72} (1950), 457--460. (Russian).



\bibitem{Pazy}
A.~Pazy,
\emph{Semigroups of Linear Operators and Applications to Partial Differential Equations},
Springer, New York 1983.

\bibitem{Pol02}
A.D.~Polyanin,
\emph{Handbook of Linear Partial Differential Equations for Engineers and Scientists},
Chapman \& Hall/CRC, 2002.

\bibitem{RS2}
M.~Reed and B.~Simon,
\emph{Methods of modern mathematical physics II: Fourier analysis, self-adjointness},
Academic Press, San Diego 1975.



\bibitem{Sims07}
R.~Sims,
\emph{Reflectionless Sturm-Liouville equations},
J. Comput. Appl.  Math.  \textbf{208} (2007), 207--225.

\bibitem{Titchmarsh}
E.C.~Titchmarsh,
\emph{Eigenfunction Expansions Associated with Second-Order Differential Equations},
Oxford University Press, London, 1946, Volumes 1 and 2.


\bibitem{Weyl}
H.~Weyl,
\emph{\"Uber gew\"ohnliche Differentialgleichungen mit Singularit\"aten und die zugeh\"origen Entwicklungen willk\"urlicher Funktionen},
Math. Ann. \textbf{68} (1910), 220--269.

\bibitem{WGK99}
M.R.M.~Witwit, N.~Gordon and J.P.~Killingbeck,
\emph{Numerical computation and analysis of the Titchmarsh-Weyl m(l) function for some simple potentials},
J. Comput. Appl. Math. \textbf{106} (1999), 131--143.

\bibitem{ZXW08}
J.~Zhang, Z.~Xu and X.~Wu,
\emph{Unified approach to split absorbing boundary conditions for nonlinear Schr\"odinger equations},
Phys. Rev. E \textbf{78} (2008), 026709.

\bibitem{ZhengFast}
C.~Zheng,
\emph{Approximation, stability and fast evaluation of exact artificial boundary condition for
the one-dimensional heat equation},
J. Comput. Math. \textbf{25} (2007), 730--745.

\bibitem{ZhengPML07}
C.~Zheng,
\emph{A perfectly matched layer approach to the nonlinear Schr\"odinger wave equations},
J.~Comput.\ Phys.\ \textbf{227} (2007), 537--556.

\end{thebibliography}
\end{document}